\documentclass[11pt,draft]{amsart}

\usepackage{a4,amssymb,amsmath,cite,epic,rotating}

\DeclareMathOperator{\Con}{Con}

\DeclareMathOperator{\Nil}{Nil}

\DeclareMathOperator{\partition}{part}

\DeclareMathOperator{\var}{var}

\newtheorem{theorem}{Theorem}[section]

\newtheorem{proposition}[theorem]{Proposition}

\newtheorem{lemma}[theorem]{Lemma}

\newtheorem{corollary}[theorem]{Corollary}

\theoremstyle{definition}

\newtheorem{remark}[theorem]{Remark}

\makeatletter

\renewcommand*\subjclass[2][2000]{\def\@subjclass{#2}\@ifundefined
{subjclassname@#1}{\ClassWarning{\@classname}{Unknown edition (#1) of
Mathematics Subject Classification; using '2000'.}}{\@xp\let\@xp
\subjclassname\csname subjclassname@#1\endcsname}}

\makeatother

\begin{document}

\title[Lower-modular elements of the lattice of semigroup varieties. III]
{Lower-modular elements\\
of the lattice of semigroup varieties. III}

\author{V. Yu. Shaprynski\v{\i}}

\address{Department of Mathematics and Mechanics, Ural State University,
Lenina 51, 620083 Ekaterinburg, Russia}

\email{vshapr@yandex.ru}

\author{B. M. Vernikov}

\address{Department of Mathematics and Mechanics, Ural State University,
Lenina 51, 620083 Ekaterinburg, Russia}

\email{boris.vernikov@usu.ru}

\date{}

\thanks{The work was partially supported by the Russian Foundation for Basic
Research (grants No.~09-01-12142,~10-01-00524) and the Federal Education
Agency of the Russian Federation (project No.~2.1.1/3537).}

\begin{abstract}
We completely determine all lower-modular elements of the lattice of all
semigroup varieties. As a corollary, we show that a lower-modular element of
this lattice is modular.
\end{abstract}



\maketitle

\section{Introduction and summary}

The collection $\mathbf{SEM}$ of all semigroup varieties forms a lattice with
respect to the class-theoretical inclusion. Special elements of different
types in this lattice have been studied in several articles. An overview of
results obtained in these articles is given in the recent survey~\cite
[Section~14]{Shevrin-Vernikov-Volkov-09}. Recall the definitions of special
elements mentioned in this paper. An element $x$ of a lattice $\langle L;
\vee,\wedge\rangle$ is called \emph{modular} if
$$\forall\,y,z\in L\colon\quad y\le z\longrightarrow(x\vee y)\wedge z=(x
\wedge z)\vee y,$$
\emph{lower-modular} if
$$\forall\,y,z\in L\colon\quad x\le y\longrightarrow x\vee(y\wedge z)=y\wedge
(x\vee z),$$
\emph{distributive} if
$$\forall y,z\in L\colon\quad x\vee(y\wedge z)=(x\vee y)\wedge(x\vee z)
\ldotp$$
\emph{Upper-modular} elements are defined dually to lower-modular ones. It is
evident that a distributive element is lower-modular.

We call a semigroup variety \emph{modular} [\emph{lower-modular,
distributive}] if it is a modular [lower-modular, distributive] element of
the lattice \textbf{SEM}. Distributive varieties are completely determined by
the authors in~\cite[Theorem~1.1]{Vernikov-Shaprynskii-distr}. Here we
consider wider class of lower-modular varieties. These varieties were
mentioned for the first time in~\cite{Vernikov-Volkov-88} (see Lemma~\ref
{0-red is cmod&lmod} below) and examined systematically in~\cite
{Vernikov-07-lmod1,Vernikov-08-lmod2}. Here we complete this examination. The
main result of this article gives a complete classification of lower-modular
varieties. To formulate this result, we need a few definitions and notation.

A pair of identities $wx=xw=w$ where the letter $x$ does not occur in the
word $w$ is usually written as the symbolic identity $w=0$. (This notation is
justified because a semigroup with such identities has a zero element and all
values of the word $w$ in this semigroup are equal to zero.) Identities of
the form $w=0$ as well as varieties given by such identities are called
0-\emph{reduced}. By $\mathcal T$, $\mathcal{SL}$, and $\mathcal{SEM}$ we
denote the trivial variety, the variety of all semilattices, and the variety
of all semigroups, respectively. The main result of the article is the
following

\begin{theorem}
\label{lmod}
A semigroup variety $\mathcal V$ is lower-modular if and only if either
$\mathcal{V=SEM}$ or $\mathcal{V=M\vee N}$ where $\mathcal M$ is one of the
varieties $\mathcal T$ or $\mathcal{SL}$, while $\mathcal N$ is a $0$-reduced
variety.
\end{theorem}

Theorem~\ref{lmod}, together with~\cite[Proposition~1.1]{Jezek-McKenzie-93}
(see also Lemmas~\ref{0-red is cmod&lmod} and~\ref{join with SL} below),
immediately implies

\begin{corollary}
\label{lmod is cmod}
A lower-modular semigroup variety is modular.\qed
\end{corollary}

Theorem~\ref{lmod} and Corollary~\ref{lmod is cmod} give answers to
Question~14.2 from~\cite{Shevrin-Vernikov-Volkov-09} and Questions~1 and~2
from~\cite{Vernikov-08-lmod2}. Besides, Theorem~\ref{lmod} solves Problems~3
and~4 from~\cite{Vernikov-08-lmod2}. It is verified in~\cite[Corollary~1.2]
{Vernikov-Shaprynskii-distr} that every distributive variety is modular.
Clearly, this claim is generalized by Corollary~\ref{lmod is cmod}.

The article consists of three sections. Section~\ref{prel} contains some
auxiliary results, while Section~\ref{proof} is devoted to the proof of
Theorem~\ref{lmod}.

\section{Preliminaries}
\label{prel}

\subsection{Some properties of lower-modular, upper-modular and modular
elements in abstract lattices and the lattice SEM}
\label{prel:lmod}

We start with two easy lattice-theoretical observations. If $L$ is a lattice
and $a\in L$ then $[a)$ stands for the \emph{principal coideal} generated by
$a$, that is, the set $\{x\in L\mid x\ge a\}$.

\begin{lemma}
\label{lmod up}
If $x$ is a lower-modular element of a lattice $L$ and $a\in L$ then the
element $x\vee a$ is a lower-modular element of the lattice $[a)$.
\end{lemma}

\begin{proof}
Let $y,z\in[a)$ and $x\vee a\le y$. Then
\begin{align*}
(x\vee a)\vee(y\wedge z)&=a\vee\bigl(x\vee(y\wedge z)\bigr)&&\\
&=a\vee\bigl(y\wedge(x\vee z)\bigr)&&\text{because $x$ is lower-modular}\\
&&&\text{and $x\le x\vee a\le y$}\\
&=y\wedge(x\vee z)&&\text{because $a\le y\wedge(x\vee z)$}\\
&=y\wedge\bigl(x\vee(a\vee z)\bigr)&&\text{because $a\le z$}\\
&=y\wedge\bigl((x\vee a)\vee z\bigr)\ldotp
\end{align*}
Thus $(x\vee a)\vee(y\wedge z)=y\wedge\bigl((x\vee a)\vee z\bigr)$, and we
are done.
\end{proof}

\begin{lemma}
\label{homomorphic image}
Let $L$ be a lattice and $\varphi$ a surjective homomorphism from $L$ onto a
lattice $L'$. If $x$ is an upper-modular element of $L$ then $\varphi(x)$ is
an upper-modular element of $L'$.
\end{lemma}

\begin{proof}
Let $x'=\varphi(x)$ and let $y',z'$ be elements of $L'$ with $y'\le x'$. Then
there are $y,z\in L$ such that $y'=\varphi(y)$ and $z'=\varphi(z)$. We may
assume that $y\le x$. Indeed, if this is not the case then we may consider
the element $x\wedge y$ rather than $y$ because $\varphi(x\wedge y)=\varphi
(x)\wedge\varphi(y)=x'\wedge y'=y'$. Since the element $x$ is upper-modular
in $L$, we have $(z\wedge x)\vee y=(z\vee y)\wedge x$. This implies that $(z'
\wedge x')\vee y'=(z'\vee y')\wedge x'$ that completes the proof.
\end{proof}

Now we provide some known partial results about lower-modular varieties. It
is well-known that if a semigroup $\mathcal V$ variety is \emph{periodic}
(that is, consists of periodic semigroups) then it contains the greatest
nil-subvariety. We denote this subvariety by $\Nil(\mathcal V)$. A semigroup
variety $\mathcal V$ is called \emph{proper} if $\mathcal{V\ne SEM}$.

\begin{lemma}[\!\!{\mdseries\cite[Theorem~1]{Vernikov-07-lmod1}}]
\label{lmod-nec}
If a proper semigroup variety $\mathcal V$ is lower-modular then $\mathcal V$
is periodic and the variety $\Nil(\mathcal V)$ is \textup0-reduced.\qed
\end{lemma}

\begin{lemma}[\!\!{\mdseries\cite[Corollary~3]{Vernikov-Volkov-88}}]
\label{0-red is cmod&lmod}
A $0$-reduced semigroup variety is modular and lower-modular.\qed
\end{lemma}

Note that the `modular half' of Lemma~\ref{0-red is cmod&lmod} was
rediscovered (in some other terms) in~\cite[Proposition~1.1]
{Jezek-McKenzie-93}.

\begin{lemma}
\label{join with SL}
A semigroup variety $\mathcal V$ is \textup[lower-\textup]modular if and only
if the variety $\mathcal{V\vee SL}$ is such.\qed
\end{lemma}

This fact was proved in~\cite[Corollary~1.5(i)]{Vernikov-Volkov-06} for
modular varieties and in~\cite[Corollary~1.3]{Vernikov-07-lmod1} for
lower-modular ones.

\subsection{Decomposition of some varieties into the join of subvarieties}
\label{prel:decomp}

We denote by $\mathcal{LZ}$ [respectively $\mathcal{RZ}$] the variety of all
left [right] zero semigroups. If $\Sigma$ is a system of semigroup identities
then $\var\Sigma$ stands for the semigroup variety given by $\Sigma$. Put
\begin{align*}
\mathcal P&=\var\{xy=x^2y,x^2y^2=y^2x^2\},\\
\overleftarrow{\mathcal P}&=\var\{xy=xy^2,x^2y^2=y^2x^2\}\ldotp
\end{align*}
Lemma~2 of the article~\cite{Volkov-89} and the proof of Proposition~1 of the
same article imply the following

\begin{lemma}
\label{M+N}
If a periodic semigroup variety $\mathcal V$ contains none of the varieties
$\mathcal{LZ}$, $\mathcal{RZ}$, $\mathcal P$, and $\overleftarrow{\mathcal
P}$ then $\mathcal{V=M\vee N}$ where the variety $\mathcal M$ is generated by
a monoid and $\mathcal{N=\Nil(V)}$.\qed
\end{lemma}

For any natural $m$, we put $\mathcal C_m=\var\{x^m=x^{m+1},xy=yx\}$. In
particular, $\mathcal C_1=\mathcal{SL}$. For notation convenience, we define
also $\mathcal C_0=\mathcal T$.

\begin{lemma}[\!\!{\mdseries\cite{Head-68}}]
\label{G+C_m}
If a semigroup variety $\mathcal M$ is generated by a commutative monoid
then $\mathcal{M=G\vee C}_m$ for some Abelian periodic group variety
$\mathcal G$ and some $m\ge0$.\qed
\end{lemma}

\subsection{Identities of certain semigroup varieties}
\label{prel:ident}

In the course of proving our results it will be convenient to have at our
disposal a description of the identities of several concrete semigroup
varieties. We denote by $F$ the free semigroup over a countably infinite
alphabet. The equality relation on $F$ is denoted by $\equiv$. If $u$ is a
word and $x$ is a letter then $c(u)$ stands for the set of all letters
occurring in $u$, $\ell(u)$ is the length of $u$, $\ell_x(u)$ denotes the
number of occurrences of $x$ in $u$, while $t(u)$ is the last letter of $u$.
The statements (i) and (ii) of the following lemma are well-known
and can be easily verified. The statement (iii) was proved in \cite[Lemma\,7]
{Golubov-Sapir-82}.

\begin{lemma}
\label{word problem}
The identity $u=v$ holds in the variety:
\begin{itemize}
\item[\textup{(i)}]$\mathcal{RZ}$ if and only if $t(u)\equiv t(v)$;
\item[\textup{(ii)}]$\mathcal C_2$ if and only if $c(u)=c(v)$ and, for every
letter $x\in c(u)$, either $\ell_x(u)>1$ and $\ell_x(v)>1$ or $\ell_x(u)=
\ell_x(v)=1$;
\item[\textup{(iii)}]$\mathcal P$ if and only if $c(u)=c(v)$ and either
$\ell_{t(u)}(u)>1$ and $\ell_{t(v)}(v)>1$ or $\ell_{t(u)}(u)=\ell_{t(v)}(v)=
1$ and $t(u)\equiv t(v)$.\qed
\end{itemize}
\end{lemma}

\subsection{Verbal subsets of free groups}
\label{prel:Sapir}

Similarly to the articles~\cite{Vernikov-07-lmod1,Vernikov-08-lmod2,
Vernikov-Shaprynskii-distr}, we need here the technique developed by Sapir
in~\cite{Sapir-91}. We introduce the basic notation from that paper. Let
$\mathcal G$ be a periodic group variety and $\{v_i=1\mid i\in I\}$ a basis
of identities of $\mathcal G$ (as a variety of groups) where $v_i$ are
semigroup words. Let $r=\exp(\mathcal G)$ where $\exp(\mathcal G)$ stands for
the exponent of the variety $\mathcal G$. For a letter $x$, put $x^0=
x^{r(r+1)}$. Let
$$S(\mathcal G)=\var\{xyz=xy^{r+1}z,\,x^0y^0=y^0x^0,\,x^2=x^{r+2},\,xv_i^2y=x
v_iy\mid i\in I\}\ldotp$$
As it is shown in~\cite{Sapir-91}, the variety $S(\mathcal G)$ does not
depend on the particular choice of the basis $\{v_i=1\mid i\in I\}$ (see
Remark~\ref{verbal subsets remark} below). Furthermore, let $F(\mathcal G)$
be the free group of countably infinite rank in $\mathcal G$. A subset $X$ of
$F(\mathcal G)$ is called \emph{verbal} if it is closed under all
endomorphisms of $F(\mathcal G)$. Clearly, a verbal subset $X$ of $F(\mathcal
G)$ is a set of all values in $F(\mathcal G)$ of some set $W$ of non-empty
words; in this case we write $X=\mathcal G(W)$. If $X$ is a verbal subset in
$F(\mathcal G)$ and $X=\mathcal G(W)$ then we put
$$S(\mathcal G,X)=S(\mathcal G)\wedge\var\bigl\{xwx=(xwx)^{r+1}\mid w\in W
\bigr\}\ldotp$$
If $X=\{1\}$ where~1 is the unit element of $F(\mathcal G)$ then we will
write $S(\mathcal G,1)$ rather than $S\bigl(\mathcal G,\{1\}\bigr)$. It is
convenient to consider the empty set as a verbal subset in $F(\mathcal G)$
and put $S(\mathcal G,\varnothing)=S(\mathcal G)$. If $\mathcal H$ is a
subvariety of $\mathcal G$ and $X$ is a verbal subset of $F(\mathcal G)$ then
we put
\begin{equation}
\label{S(H,X)}
S(\mathcal H,X)=S(\mathcal H)\wedge S(\mathcal G,X)\ldotp
\end{equation}
To avoid a possible confusion, we note that the paper~\cite{Sapir-91} does
not contain an explicit definition of the variety $S(\mathcal H,X)$ where $X$
is a verbal subset of $F(\mathcal G)$ with $\mathcal{G\ne X}$. But one can
trace the argument of~\cite{Sapir-91} to see that the equality~\eqref{S(H,X)}
is what Sapir tacitly meant by this definition but failed to explicitly
define.

As usual, if $\mathcal X$ is a variety then $L(\mathcal X)$ stands for the
subvariety lattice of $\mathcal X$. To prove Theorem~\ref{lmod}, we need the
following

\begin{lemma}[\!\!{\mdseries\cite{Sapir-91}}]
\label{verbal subsets}
Let $\mathcal G$ be a variety of periodic groups. The interval $\bigl[S
(\mathcal T,1)$, $S(\mathcal G)\bigr]$ of the lattice $L\bigl(S(\mathcal G)
\bigr)$ consists of all varieties of the form $S(\mathcal H,X)$ where
$\mathcal{H\subseteq G}$ and $X$ is a \textup(possibly empty\textup) verbal
subset of $F(\mathcal G)$. Here, for varieties $S(\mathcal H,X)$ and $S
(\mathcal H',X')$ from the interval $\bigl[S(\mathcal T,1),\,S(\mathcal G)
\bigr]$, the inclusion $S(\mathcal H',X')\subseteq S(\mathcal H,X)$ holds if
and only if $\mathcal{H'\subseteq H}$ and there exists a set of words $W$
such that $X=\mathcal H(W)$ and $\mathcal H'(W)\subseteq X'$.\qed
\end{lemma}

\begin{remark}
\label{verbal subsets remark}
Lemma~\ref{verbal subsets} shows that the construction of the variety $S
(\mathcal G,X)$ is in fact independent of the actual choice of the
`generator' $W$ of the verbal subset $X$; it is only $X$ that really matters,
as different choices of $W$ will result in the same variety. In particular,
by the definition of the variety $S(\mathcal G,X)$, it satisfies the identity
$xwx=(xwx)^{r+1}$ whenever $w\in W$. In view of Lemma~\ref{verbal subsets},
this identity holds in $S(\mathcal G,X)$ not only for $w\in W$ but for any
word $w$ representing an element of $X$.
\end{remark}

\subsection{Overcommutative varieties}
\label{prel:oc}

We denote by $\mathcal{COM}$ the variety of all commutative semigroups. A
semigroup variety $\mathcal V$ is called \emph{overcommutative} if $\mathcal
{V\supseteq COM}$. The lattice of all overcommutative varieties is denoted by
\textbf{OC}. The structure of this lattice was clarified by Volkov in~\cite
{Volkov-94}. It turns out that the lattice $\mathbf{OC}$ admits a concise and
transparent description in terms of congruence lattices of unary algebras of
some special type, called $G$-sets. This description plays an essential role
in the proof of Theorem~\ref{lmod}. To reproduce the result from~\cite
{Volkov-94}, we need some new definitions and notation.

Let $A$ be a non-empty set. We denote by $\mathbf S_A$ the group of all
permutations on $A$. If $A=\{1,2,\dots,m\}$ then we will write $\mathbf S_m$
rather than $\mathbf S_{\{1,2,\dots,m\}}$. A $G$-\emph{set} is a unary
algebra on a set $A$ where the unary operations form a group of permutations
on $A$ (that is, the subgroup of $\mathbf S_A$). The congruence lattice of a
$G$-set $A$ is denoted by $\Con(A)$.

Let $m$ and $n$ be positive integers with $2\le m\le n$. A sequence $\lambda=
(\ell_1,\ell_2,\dots,\ell_m)$ of positive integers such that
$$\ell_1\ge\ell_2\ge\cdots\ge\ell_m\enskip\text{and}\enskip\sum_{i=1}^m
\ell_i=n$$
is said to be a \emph{partition of the number $n$ into $m$ parts}. For a word
$u$, we put $\partition(u)=\bigl(\ell_{x_1}(u),\ell_{x_2}(u),\dots,\ell_{x_m}
(u)\bigr)$ where $m=\max\bigl\{i\mid x_i\in c(u)\bigr\}$. Let us fix positive
integers $m$ and $n$ with $2\le m\le n$ and a partition $\lambda=(\ell_1,
\ell_2,\dots,\ell_m)$ of the number $n$ into $m$ parts. Put
\begin{align*}
W_\lambda&=\bigl\{u\in F\mid\ell(u)=n,\,c(u)=\{x_1,x_2,\dots,x_m\}\ \text
{and}\ \partition(u)=\lambda\bigr\},\\
\mathbf S_\lambda&=\{\sigma\in\mathbf S_m\mid\ell_i=\ell_{i\sigma}\ \text
{for all}\ i=1,2,\dots,m\}\ldotp
\end{align*}
Clearly, $\mathbf S_\lambda$ is a subgroup in $\mathbf S_m$.

If $u\equiv x_{i_1}x_{i_2}\cdots x_{i_n}$ where $x_{i_1},x_{i_2},\dots,
x_{i_n}$ are (not necessarily different) letters and $\pi\in\mathbf S_{c(u)}$
then we denote by $u\pi$ the word $x_{i_1\pi}x_{i_2\pi}\cdots x_{i_n\pi}$. It
is clear that if $u\in W_\lambda$ and $\pi\in\mathbf S_\lambda$ then $u\pi\in
W_\lambda$. For every $\pi\in\mathbf S_\lambda$, we define the unary
operation $\pi^*$ on $W_\lambda$ by letting $\pi^*(u)\equiv u\pi$ for any
word $u\in W_\lambda$. Obviously, the set $W_\lambda$ with the collection of
unary operations $\{\pi^*\mid\pi\in\mathbf S_\lambda\}$ is an $\mathbf
S_\lambda$-set. The description of the lattice \textbf{OC} mentioned above is
given by the following

\begin{proposition}[\!\!{\mdseries\cite{Volkov-94}}]
\label{oc:struct}
The lattice $\mathbf{OC}$ is anti-isomorphic to a subdirect product of
congruence lattices $\Con(W_\lambda)$ where $\lambda$ runs over the set of
all partitions.\qed
\end{proposition}

\section{Proof of Theorem \ref{lmod}}
\label{proof}

\emph{Sufficiency} immediately follows from Lemmas~\ref{0-red is cmod&lmod}
and~\ref{join with SL} and the evident fact that the variety $\mathcal{SEM}$
is lower-modular.

\smallskip

\emph{Necessity}. Let $\mathcal V$ be a proper lower-modular semigroup
variety. Lemma~\ref{lmod up} implies that the variety $\mathcal{W=V\vee COM}$
is a lower-modular element of the lattice \textbf{OC}. The variety $\mathcal
W$ is proper because the variety $\mathcal{SEM}$ is not decomposable into the
join of any two proper varieties~\cite{Dean-Evans-69}.

Recall that an identity $u=v$ is called \emph{balanced} if each letter occurs
in $u$ and $v$ the same number of times. It is well-known that if an
overcommutative variety satisfies some identity then this identity is
balanced.

Being proper, the variety $\mathcal W$ satisfies a non-trivial balanced
identity $u=v$. Let $|c(u)|=m$. We may assume that $c(u)=\{x_1,x_2,\dots,
x_m\}$ and $\ell_{x_1}(u)\ge\ell_{x_2}(u)\ge\cdots\ge\ell_{x_m}(u)$
(otherwise we may rename letters). Put $\ell_i=\ell_{x_i}(u)$ for all $i=1,2,
\dots,m$. Then $\partition(u)=\partition(v)=(\ell_1,\ell_2,\dots,\ell_m)$. We
may assume that $\ell_1>\ell_2>\cdots>\ell_m>1$ (if it is not the case, we
may multiply $u=v$ by an appropriate word on the right). Let $x$ and $y$ be
arbitrary letters with $x,y\notin c(u)$ and $\lambda=\partition(xyu)=(\ell_1,
\ell_2,\dots,\ell_m,1,1)$ (we identify here $x$ and $y$ with $x_{m+1}$ and
$x_{m+2}$ respectively). We denote by $\nu$ the fully invariant congruence on
$F$ corresponding to the variety $\mathcal W$ and by $\alpha$ the restriction
of $\nu$ to $W_\lambda$. Then
$$xyu\,\alpha\,xyv,\,yxu\,\alpha\,yxv,\,xuy\,\alpha\,xvy,\,yux\,\alpha\,yvx
\ldotp$$
Proposition~\ref{oc:struct} implies that there is a surjective homomorohism
from the lattice dual to \textbf{OC} onto $\Con(W_\lambda)$. Now Lemma~\ref
{homomorphic image} applies with the conclusion that $\alpha$ is an
upper-modular element of the lattice $\Con(W_\lambda)$. Since $\ell_1>\ell_2>
\cdots>\ell_m>1$, the group $\mathbf S_\lambda$ consists of two elements. Let
$\beta$ be the equivalence relation on $W_\lambda$ with only two
non-singleton classes $\{xyu,xyv\}$ and $\{yxu,yxv\}$, and let $\gamma$ be
the equivalence relation on $W_\lambda$ with only four non-singleton classes
$\{xyu,xuy\}$, $\{xyv,xvy\}$, $\{yxu,yux\}$, and $\{yxv,yvx\}$. It is evident
that $\beta$ and $\gamma$ are congruences on $W_\lambda$ and $\beta\subseteq
\alpha$. Therefore $(\gamma\vee\beta)\wedge\alpha=(\gamma\wedge\alpha)\vee
\beta$.

Note that $xuy\,\gamma\,xyu\,\beta\,xyv\,\gamma\,xvy$. Therefore $(xuy,xvy)
\in\gamma\vee\beta$, whence
$$(xuy,xvy)\in(\gamma\vee\beta)\wedge\alpha=(\gamma\wedge\alpha)\vee\beta
\ldotp$$
This means that there is a word $w\in W_\lambda$ such that $xuy\not\equiv w$
and either $(xuy,w)\in\gamma\wedge\alpha$ or $xuy\,\beta\,w$. But the latter
contradicts the choice of the congruence $\beta$. Hence $(xuy,w)\in\gamma
\wedge\alpha$. In particular, $xuy\,\gamma\,w$. The definition of $\gamma$
implies that $w\equiv xyu$. Thus $(xuy,xyu)\in\gamma\wedge\alpha$. In
particular, we have verified that $xuy\,\alpha\,xyu$.

Let $\gamma'$ be the equivalence relation on $W_\lambda$ with only four
non-singleton classes $\{xyu,yux\}$, $\{xyv,yvx\}$, $\{yxu,xuy\}$, and $\{yx
v,xvy\}$. Clearly, $\gamma'$ is a congruence on $W_\lambda$. Now we may
repeat almost literally arguments from the previous paragraph with using
$\gamma'$ rather than $\gamma$. As a result, we obtain that $xyu\,\alpha\,yu
x$. Thus $xuy\,\alpha\,xyu\,\alpha\,yux$. Since $\alpha\subseteq\nu$, we have
$xuy\,\nu\,yux$. This means that the identity
\begin{equation}
\label{xuy=yux}
xuy=yux
\end{equation}
holds in the variety $\mathcal W$. Therefore this identity holds in the
variety $\mathcal V$ as well.

Lemma~\ref{word problem} and its dual imply that the identity~\eqref{xuy=yux}
fails in the varieties $\mathcal{LZ}$, $\mathcal{RZ}$, $\mathcal P$, and
$\overleftarrow{\mathcal P}$. Thus $\mathcal V$ does not contain these
varieties. By Lemma~\ref{lmod-nec} the variety $\mathcal V$ is periodic. Now
Lemma~\ref{M+N} applies with the conclusion that $\mathcal{V=M\vee N}$ where
the variety $\mathcal M$ is generated by a monoid and $\mathcal{N=\Nil(V)}$.
Lemma~\ref{lmod-nec} implies that the variety $\mathcal N$ is 0-reduced. It
remains to verify that $\mathcal M$ is one of the varieties $\mathcal T$ or
$\mathcal{SL}$.

Since $\mathcal M$ is generated by a monoid, the set of its identities is
closed for deleting letters; therefore, by deleting of all letters from $c
(u)$ in~\eqref{xuy=yux}, we obtain that $\mathcal M$ is commutative. Now we
can apply Lemma~\ref{G+C_m} and conclude that $\mathcal{M=G\vee C}_m$ for
some Abelian periodic group variety $\mathcal G$ and some $m\ge0$. Suppose
that $m\ge2$. Then $\mathcal{V\supseteq C}_2$. It is easy to deduce from
Lemma~\ref{word problem} that $\mathcal C_2\vee\mathcal{RZ\supseteq P}$.
Hence $\mathcal{V\vee RZ\supseteq C}_2\vee\mathcal{RZ\supseteq P}$. Put
$\mathcal{U=V\vee P}$. Note that $\mathcal{V,P\nsupseteq RZ}$. As is
well-known, the variety $\mathcal{RZ}$ is an atom of the lattice \textbf
{SEM}. Therefore $\mathcal{V\wedge RZ=P\wedge RZ=T}$. It is well-known also
that the lattice \textbf{SEM} is 0-\emph{distributive}, that is, satisfies
the condition
$$\forall x,y,z\colon\quad x\wedge z=y\wedge z=0\longrightarrow(x\vee y)
\wedge z=0$$
(see~\cite[Section~1]{Shevrin-Vernikov-Volkov-09}, for instance). Therefore
$\mathcal{U\wedge RZ=(V\vee P)\wedge RZ=T}$. Combining these observations, we
have
\begin{align*}
\mathcal V&=\mathcal{V\vee T}&&\\
&=\mathcal{V\vee(U\wedge RZ)}&&\text{because $\mathcal{U\wedge RZ=T}$}\\
&=\mathcal{U\wedge(V\vee RZ)}&&\text{because $\mathcal V$ is lower-modular
and $\mathcal{V\subseteq U}$}\\
&\supseteq\mathcal P&&\text{because $\mathcal{U=V\vee P\supseteq P}$ and
$\mathcal{V\vee RZ\supseteq P}$}\ldotp
\end{align*}
Thus $\mathcal{V\supseteq P}$. A contradiction shows that $m\le1$, whence
$\mathcal{M=G\vee K}$ where $\mathcal K$ is one of the varieties $\mathcal T$
or $\mathcal{SL}$. It remains to check that $\mathcal{G=T}$.

The remaining part of the proof is based on Lemma~\ref{verbal subsets}. Note
that in what follows we combine (with slight modifications) arguments from
the proofs of~\cite[Lemma~2.2]{Vernikov-08-lmod2} and~\cite[Proposition~3.1]
{Vernikov-Shaprynskii-distr}.

\begin{figure}[tbh]
\begin{center}
\unitlength=1mm
\linethickness{0.4pt}
\begin{picture}(42,56)
\put(1,23){\line(1,-1){20}}
\put(1,23){\line(1,1){10}}
\put(11,13){\line(1,1){10}}
\put(11,33){\line(0,1){10}}
\put(11,33){\line(1,-1){10}}
\put(11,43){\line(1,1){10}}
\put(11,43){\line(1,-1){10}}
\put(21,3){\line(1,1){20}}
\put(21,23){\line(0,1){10}}
\put(21,23){\line(1,-1){10}}
\put(21,33){\line(1,1){10}}
\put(21,53){\line(1,-1){10}}
\put(31,43){\line(1,-2){10}}
\put(1,23){\circle*{1.33}}
\put(11,13){\circle*{1.33}}
\put(11,33){\circle*{1.33}}
\put(11,43){\circle*{1.33}}
\put(21,3){\circle*{1.33}}
\put(21,23){\circle*{1.33}}
\put(21,33){\circle*{1.33}}
\put(21,53){\circle*{1.33}}
\put(31,13){\circle*{1.33}}
\put(31,43){\circle*{1.33}}
\put(41,23){\circle*{1.33}}
\put(9,10){\makebox(0,0)[cc]{$\mathcal G$}}
\put(33,46){\makebox(0,0)[cc]{$S(\mathcal G)$}}
\put(21,36){\makebox(0,0)[cc]{$S(\mathcal G,1)$}}
\put(21,26){\makebox(0,0)[cc]{$S\bigl(\mathcal G,F(\mathcal G)\bigr)$}}
\put(42,23){\makebox(0,0)[lc]{$S(\mathcal T)$}}
\put(36,10){\makebox(0,0)[cc]{$S(\mathcal T,1)$}}
\put(0,23){\makebox(0,0)[rc]{$\mathcal V$}}
\put(21,56){\makebox(0,0)[cc]{$\mathcal V\vee S(\mathcal G)$}}
\put(10,43){\makebox(0,0)[rc]{$\mathcal V\vee S(\mathcal G,1)$}}
\put(10,33){\makebox(0,0)[rc]{$\mathcal V\vee S\bigl(\mathcal G,F(\mathcal G)
\bigr)$}}
\put(21,0){\makebox(0,0)[cc]{$\mathcal T$}}
\end{picture}
\caption{A fragment of the lattice $L\bigl(\mathcal V\vee S(\mathcal G)
\bigr)$}
\label{V+S(G)}
\end{center}
\end{figure}
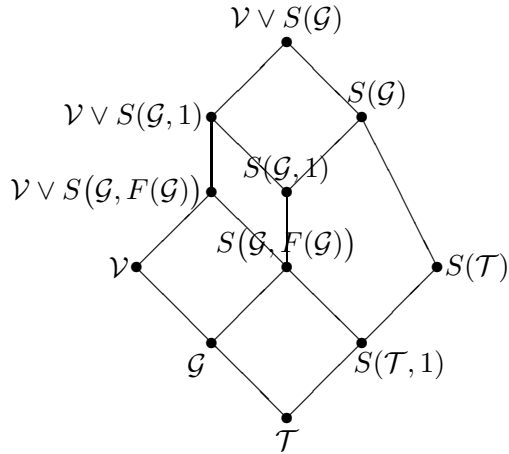

Suppose that $\mathcal{G\ne T}$. Put $\mathcal{Y=V}\vee S(\mathcal G,1)$ and
$\mathcal Z=S(\mathcal T)$. Further considerations are illustrated by Fig.~\ref
{V+S(G)}.

Lemma~\ref{verbal subsets} implies that $S(\mathcal G)=S\mathcal{(T)\vee G}$
(see Fig.~\ref{V+S(G)}). Using this equality and the inclusion $\mathcal{G
\subseteq V}$, we have
$$\mathcal Y=S(\mathcal G,1)\vee\mathcal V\subseteq S\mathcal{(G)\vee V}=S
\mathcal{(T)\vee G\vee V}=S\mathcal{(T)\vee V=Z\vee V}\ldotp$$
Therefore $\mathcal{(Z\vee V)\wedge Y=Y}$. Since the variety $\mathcal V$ is
lower-modular and $\mathcal{V\subseteq Y}$, we have $\mathcal{(Z\wedge Y)\vee
V=(Z\vee V)\wedge Y}$, whence
\begin{equation}
\label{ZY+V=Y}
\mathcal{(Z\wedge Y)\vee V=Y}\ldotp
\end{equation}
Furthermore, Lemma~\ref{verbal subsets} implies that $S(\mathcal T,1)=S
(\mathcal T)\wedge S(\mathcal G,1)$ (see Fig.~\ref{V+S(G)}). Therefore
$$S(\mathcal T,1)=S(\mathcal T)\wedge S(\mathcal G,1)\subseteq S\mathcal
{\bigl(T)\wedge(V}\vee S(\mathcal G,1)\bigr)=\mathcal{Z\wedge Y\subseteq Z}=S
(\mathcal T),$$
that is, $S(\mathcal T,1)\subseteq\mathcal{Z\wedge Y}\subseteq S(\mathcal
T)$. It is evident that the group $F(\mathcal T)$ contains only two verbal
subsets, namely $\varnothing$ and \{1\}. Therefore Lemma~\ref{verbal subsets}
implies that the interval $\bigl[S(\mathcal T,1),\,S(\mathcal T)\bigr]$ of
the lattice $L\bigl(S(\mathcal T)\bigr)$ consists of the varieties $S
(\mathcal T,1)$ and $S(\mathcal T)$ only. Thus either $\mathcal{Z\wedge Y}=S
(\mathcal T,1)$ or $\mathcal{Z\wedge Y}=S(\mathcal T)$. Let us consider these
two cases separately. Let $\exp(\mathcal G)=r$.

\smallskip

\emph{Case}~1: $\mathcal{Z\wedge Y}=S(\mathcal T,1)$. Lemma~\ref
{verbal subsets} implies that $S(\mathcal T,1)\vee\mathcal G=S\bigl(\mathcal
G,F(\mathcal G)\bigr)$ (see Fig.~\ref{V+S(G)}). Using the equality~\eqref
{ZY+V=Y} and the inclusion $\mathcal{G\subseteq V}$, we have
\begin{align*}
S\bigl(\mathcal G,F\mathcal{(G)\bigr)\vee V}&=S(\mathcal T,1)\vee\mathcal{G
\vee V}=S(\mathcal T,1)\vee\mathcal V\\
&=\mathcal{(Z\wedge Y)\vee V=Y}=S(\mathcal G,1)\vee\mathcal V,
\end{align*}
that is,
\begin{equation}
\label{S(G,F(G))+V=S(G,1)+V}
S\bigl(\mathcal G,F\mathcal{(G)\bigr)\vee V}=S(\mathcal G,1)\vee\mathcal V
\ldotp
\end{equation}
Being a nilvariety, $\mathcal N$ satisfies an identity $u=0$ for some word
$u$. Suppose that the variety $\mathcal G$ (considered as a variety of
groups) satisfies the identity $u=1$. Let $x$ be a letter with $x\notin c
(u)$. Since the variety $\mathcal G$ is non-trivial, it does not satisfy the
identity $ux=1$. It is evident that $ux=0$ in $\mathcal N$. Thus there is a
word $w$ such that the variety $\mathcal N$ satisfies the identity $w=0$ but
the variety $\mathcal G$ does not satisfy the identity $w=1$. Let $x$ be a
letter with $x\notin c(w)$. Remark~\ref{verbal subsets remark} implies that
$S\bigl(\mathcal G,F(\mathcal G)\bigr)$ satisfies the identity
\begin{equation}
\label{key ident1}
xwx=(xwx)^{r+1}\ldotp
\end{equation}
This identity holds in the variety $\mathcal V$ as well because $\mathcal{V
\subseteq G\vee SL\vee N}$. Therefore the variety $\mathcal V\vee S\bigl
(\mathcal G,F(\mathcal G)\bigr)$ satisfies the identity~\eqref{key ident1}.
But~\eqref{key ident1} fails in the variety $S(\mathcal G,1)$ by the
definition of this variety, whence~\eqref{key ident1} fails in the variety
$\mathcal V\vee S(\mathcal G,1)$. We have a contradiction with the
equality~\eqref{S(G,F(G))+V=S(G,1)+V}.

\smallskip

\emph{Case}~2: $\mathcal{Z\wedge Y}=S(\mathcal T)$. As we have already noted
above, $S(\mathcal G)=S\mathcal{(T)\vee G}$ (see Fig.~\ref{V+S(G)}). Taking
into account the equality~\eqref{ZY+V=Y} and the inclusion $\mathcal{G
\subseteq V}$, we have
$$S(\mathcal G,1)\vee\mathcal{V=Y=(Z\wedge Y)\vee V}=S\mathcal{(T)\vee V}=
S\mathcal{(T)\vee G\vee V}=S\mathcal{(G)\vee V}\ldotp$$
We see that
\begin{equation}
\label{S(G,1)+V=S(G)+V}
S(\mathcal G,1)\vee\mathcal V=S(\mathcal G)\vee\mathcal V\ldotp
\end{equation}
Let $w$ be an arbitrary word such that the variety $\mathcal G$ satisfies (as
a variety of groups) the identity $w=1$. Being a nil-variety, $\mathcal N$
satisfies the identity $x^n=0$ for some $n$. The variety $\mathcal G$
satisfies the identity $w^{2n}=1$. Remark~\ref{verbal subsets remark} implies
that the variety $S(\mathcal G,1)$ satisfies the identity
\begin{equation}
\label{key ident2}
xw^{2n}x=(xw^{2n}x)^{r+1}\ldotp
\end{equation}
This identity holds in the varieties $\mathcal M$ and $\mathcal N$ as well,
whence it holds in $\mathcal V$, and therefore in $S(\mathcal G,1)\vee
\mathcal V$. The equality~\eqref{S(G,1)+V=S(G)+V} implies that~\eqref
{key ident2} holds in $S\mathcal{(G)\vee V}$, and therefore in $S(\mathcal
G)$. We always may include the identity $w=1$ in the identity basis of
$\mathcal G$. By the definition of the variety $S(\mathcal G)$, it satisfies
the identity $xwx=xw^2x$, and therefore the identities
\begin{align*}
xwx&=xw^2x=xw^4x\equiv x(w\cdot w^2\cdot w)x=x(w\cdot w^4\cdot w)x\equiv xw^6
x\\
&\equiv x(w^2\cdot w^2\cdot w^2)x=x(w^2\cdot w^4\cdot w^2)x\equiv
xw^8x=\cdots=xw^{2n}x\ldotp
\end{align*}
Combining the identities $xwx=xw^{2n}x$ and~\eqref{key ident2}, we have that
the identities
$$xwx=xw^{2n}x=(xw^{2n}x)^{r+1}=(xwx)^{r+1}$$
hold in $S(\mathcal G)$. Thus $S(\mathcal G)$ satisfies the identity~\eqref
{key ident1} whenever $\mathcal G$ satisfies $w=1$. Therefore $S(\mathcal G)
\subseteq S(\mathcal G,1)$ but this inclusion contradicts Lemma~\ref
{verbal subsets}.

\smallskip

We have verified that $\mathcal{G=T}$ and completed the proof of Theorem~\ref
{lmod}.\qed

\subsection*{Acknowledgements}
The authors would like to express their gratitude to the anonymous referee
for valuable remarks and suggestions on improving the paper.

\end{document}